\documentclass[12pt,dvips]{article}
\usepackage{amsmath,amsfonts,amssymb,amsthm,graphicx,verbatim,subfig}
\usepackage[all]{xy}

\ifx\pdftexversion\undefined

\usepackage[a4paper,colorlinks,link
=black,filecolor=black,citecolor=black,urlcolor=black,pdfstartview=FitH]{hyperref}
\else

\usepackage[a4paper,colorlinks,linkcolor=black,filecolor=black,citecolor=black,urlcolor=black,pdfstartview=FitH]{hyperref}
\fi

\usepackage[all]{xy}
\usepackage[a4paper,left=25mm]{geometry}
\ifx\pdftexversion\undefined

\usepackage[a4paper,colorlinks,linkcolor=black,filecolor=black,citecolor=black,urlcolor=black,pdfstartview=FitH]{hyperref}
\fi

\font\sixbb=msbm6
\font\eightbb=msbm8
\font\twelvebb=msbm10 scaled 1095
\newfam\bbfam
\textfont\bbfam=\twelvebb \scriptfont\bbfam=\eightbb
                           \scriptscriptfont\bbfam=\sixbb

\newcommand{\FF}{\mathbb{F}}

\newtheorem{theorem}{\bf Theorem}[section]

\newtheorem{proposition}[theorem]{\bf Proposition}

\newcommand{\enp}{\begin{flushright} $\Box$ \end{flushright}}
\newcommand{\beq}[0]{\begin{equation}}
\newcommand{\enq}[0]{\end{equation}}

\newcommand{\supp}{{\rm supp}}
\newcommand{\cf}{{\cal F}}

\newcommand{\dist}{{\rm dist}}
\newcommand{\calc}{\mathcal{C}}

\title{Graph Codes and Local Systems}
\begin{document}
\author{Roy Meshulam\thanks{Department of Mathematics,
Technion, Haifa 32000, Israel. e-mail:
meshulam@math.technion.ac.il~. Supported by ISF and GIF grants.}}
\maketitle
\pagestyle{plain}
\begin{abstract}
It is shown that the good expander codes introduced by Sipser and Spielman, can be realized as the first homology of a graph with respect to a certain twisted coefficient system.
\end{abstract}

\section{Introduction}
\label{s:intro}

We first recall some standard definitions from coding theory.
The {\it Hamming weight} of $w=(w_1,\ldots,w_N) \in \FF_2^N$ is
$|w|:=|\{1 \leq i \leq N: w_i \neq 0\}|$.
A {\it binary linear code} of length $N$ is a linear subspace
$C \subset \FF_2^N$. The {\it distance} of $C$ is
$$\dist(C):= \min\{|w|:0 \neq w \in C\}.$$
The {\it rate} and {\it relative distance} of $C$ are given respectively by
$$r(C):=\frac{\dim C}{N}~~~~,~~~~ \delta(C):=\frac{\dist(C)}{N}.$$
Let $\{N_i\}_{i=1}^{\infty}$ be a sequence of integers
such that $N_i \rightarrow \infty$. A family of binary linear codes $C_i \subset \FF_2^{N_i}$ is
{\it asymptotically good} if
there exist $r, \delta>0$ such that $r(C_i) \geq r$ and $\delta(C_i)\geq \delta$ for all $i$.
A classical probabilistic argument due to Shannon and to Gilbert and Varshamov shows that there exist families of asymptotically good linear codes. Explicit constructions of such families are more involved and were first given by Justesen.
The study of good codes remains a major area of research in coding theory to this day. See van Lint's book \cite{vanLint99} for a comprehensive introduction to the subject.

Here we are concerned with codes that arise from simplicial complexes.
Let $X$ be an $n$-dimensional simplicial complex with $N=f_n(X)$ top dimensional faces.
Let $\sigma_1,\ldots,\sigma_N$ denote the $n$-simplices of $X$.
We identify the space of $n$-chains $C_n(X,\FF_2)$ with $\FF_2^N$,
via the mapping $\sum_{i=1}^N a_i \sigma_i \rightarrow (a_1,\ldots,a_N)$.
The {\it homological code} associated to $X$ is the
homology space $C=H_n(X;\FF_2) \subset \FF_2^N$.

For $n=1$, the homological code associated to a graph $G$ is just the cycle space of $G$, and
it is well known that this code is not good. In fact, a strong form of the Moore bound due to Alon, Hoory and Linial \cite{AHL02} implies that for fixed $\epsilon>0$, if a graph $X$ on $N$ vertices has $(1+\epsilon)N$ edges, then $X$ contains a cycle of size $O_{\epsilon}(\log N)$.
The situation is different for $n \geq 2$, where the existence of good homological codes associated
with $n$-dimensional complexes follows from Theorem 4.1 in \cite{ALLM13}
(see also Theorem 6.3(1) in  \cite{DGK18}). Let $\Delta_{N-1}^{(n)}$ denote the $n$-dimensional skeleton of the $(N-1)$-simplex $\Delta_{N-1}$.
\begin{theorem}[\cite{ALLM13}]
\label{t:goodc}
Let $n \geq 2$ be fixed. Then for any integer $K >0$ there is a $\lambda=\lambda(n,K)>0$,
such that for sufficiently large $N$ there exists a complex $X_N \subset \Delta_{N-1}^{(n)}$
such that $f_n(X_N)= K \binom{N}{n}$ and $|\supp(z)| \geq \lambda \binom{N}{n}$ for all
$0 \neq z \in C=H_n(X_N;\FF_2)$.
In particular, $r(C) \geq \frac{K-1}{K}$ and $\delta(C) \geq \frac{\lambda}{K}$.
\end{theorem}
\noindent
Returning to the $1$-dimensional case, while the cycle space of a graph is not a good code,
there is a remarkable construction of graph based codes, due to Tanner \cite{Tanner81} and Sipser and Spielman \cite{SiSp96}, that does give rise to good codes. These codes are defined as follows. Let $G=(V,E)$ be a graph. For $u \in V$ let $\Gamma_G(u)=\{v \in V: uv \in E\}$. Suppose that for each vertex $u \in V$ there is given a linear subspace $C_u \subset \FF_2^{\Gamma_G(u)}$ (the local code). The {\it graph code}
$\calc=\calc(G, \{C_u\}_{u \in V})$
 is the linear space
$$\calc=\{x=(x_{uv})_{uv \in E} \in \FF_2^E: (x_{uv}: v \in \Gamma_G(u)) \in C_u {\rm~~for~all~~} u \in V\}.$$
Let $\lambda(G)$ denote the second largest eigenvalue of the adjacency matrix of $G$.
\begin{theorem}[Sipser and Spielman \cite{SiSp96}]
\label{t:sisp}
Suppose $G$ is a $d$-regular graph with $\lambda(G) \leq \lambda$, and suppose the codes $C_u$ satisfy $r(C_u) \geq r > \frac{1}{2}$ and $\delta(C_u) \geq \delta$ for all $u \in V$.
Then the graph code $\calc$ satisfies $r(\calc) \geq 2r-1$ and
$$\delta(\calc) \geq \left(\frac{\delta-\frac{\lambda}{d}}{1-\frac{\lambda}{d}}\right)^2.$$
\end{theorem}
\noindent
Let $d \geq 3$ and let $\{G_n=(V_n,E_n)\}_n$ be a sequence of $d$-regular expanders on $n$ vertices such that
$\frac{\lambda(G_n)}{d} \rightarrow 0$. For each $u \in V_n$ choose a linear subspace $C_{n,u} \subset \FF_2^{\Gamma_{G_n}(u)}$
such that $\dim C_{n,u} > \frac{r}{2}$. Theorem \ref{t:sisp} implies that
$\{\calc(G_n, \{C_{n,u}\}_{u \in V_n})\}_n$ is a family of asymptotically good codes.
For further salient properties of these codes (e.g. linear decoding time) see \cite{SiSp96,Sp99}.

In view of the importance of the Sipser-Spielman codes, it might be worthwhile
to observe that a graph code $\calc=\calc(G, \{C_{u}\}_{u \in V})$
can be obtained as the first homology $H_1(G;\cf)$ of the graph $G$ with coefficients in certain local system $\cf$ that depends on the local codes $\{C_{u}\}_{u \in V}$. In the following section we recall the relevant homological notions and describe the construction of $\cf$.

\section{Graph Codes via Twisted Homology}
\label{s:twist}

Let $X$ be a finite simplicial complex on the vertex set $V$. Let $\prec$ be an arbitrary but fixed
linear order on $V$. Let $X(k)$ denote the set of $k$-dimensional simplicies of $X$.
For a $\sigma=\{v_0,\ldots,v_k\} \in X(k)$ where $v_0 \prec \cdots \prec v_k$, we will write
$\sigma=[v_0 \ldots v_k]$ or, when there is no danger of confusion, $\sigma=v_0\ldots v_{k+1}$.
The $i$-the face of $\sigma$ as above is defined by
$\sigma_i=\sigma \setminus \{v_i\}$.
A {\it local system} $\cf$ on $X$ (see e.g. section 7 in \cite{Quillen78}) is an assignment of an abelian group $\cf(\sigma)$ to each simplex $\sigma \in X$, together with homomorphisms $\rho_{\sigma}^{\tau}:\cf(\tau) \rightarrow \cf(\sigma)$ for each
$\sigma \subset \tau$ satisfying the usual compatibility conditions: $\rho_{\sigma}^{\sigma}={\rm identity}$, and $\rho_{\eta}^{\sigma} \rho_{\sigma}^{\tau}=\rho_{\eta}^{\tau}$ if $\eta \subset \sigma \subset \tau$.
Let $$C_k(X;\cf)=\bigoplus_{\sigma \in X(k)} \cf(\sigma).$$
Elements of $C_k(X;\cf)$ are written as $c=\sum_{\sigma \in X(k)} a_{\sigma} \sigma$ with $a_{\sigma} \in \cf(\sigma)$. Define the boundary map $$\partial_k:C_k(X;\cf) \rightarrow C_{k-1}(X;\cf)$$ by
$$\partial_k c=\sum_{\sigma \in X(k)} \sum_{i=0}^k (-1)^i \rho_{\sigma_i}^{\sigma}(a_{\sigma}) \sigma_i.$$
The homology of the complex $C_*(X;\cf)$ is denoted by $H_*(X;\cf)$.

Let $G=(V,E)$ be a graph and for $u \in V$ let $C_u \subset \FF_2^{\Gamma_G(u)}$ be the local code associated to $u$. We now show that the graph code $\calc=\calc(G, \{C_u\}_{u \in V})$ can be realized as $H_1(G;\cf)$ for a certain local system $\cf$.
Let $A_u$ be a parity check matrix of $C_u$, and
let $\{w_{uv}:v \in \Gamma_G(u)\}$ be the columns of $A_u$. Thus, a vector
$y=(y_v)_{v \in \Gamma_G(u)} \in \FF_2^{\Gamma_G(u)}$
satisfies $y \in C_u$ iff $\sum_{v \in \Gamma_G(u)} y_v w_{uv} =0$.
Let $L_u=\{A_u z: z \in \FF_2^{\Gamma_G(u)} \}$ denote the column space of $A_u$.
Define a local system $\cf$ on $G$ as follows:
$$
\cf(\sigma)=
\left\{
\begin{array}{ll}
L_u & \sigma=u \in V, \\
\FF_2 & \sigma=uv \in E.
\end{array}
\right.~~
$$
The map $\rho_{u}^{uv}:\cf(uv)=\FF_2 \rightarrow \cf(u)=L_u$ is given by
$\rho_u^{uv}(\epsilon)=\epsilon w_{uv}$.
\begin{proposition}
\label{p:clocal}
$$\calc=H_1(G;\cf).$$
\end{proposition}
\noindent
{\bf Proof.} Let $c= \sum_{uv \in E} x_{uv} [uv] \in C_1(G;\cf)$.
As our ground field is $\FF_2$, we have:
\begin{equation*}
\label{e:part1}
\begin{split}
\partial_1 c &=\partial_1 \left( \sum_{uv \in E} x_{uv} [uv] \right) \\
&= \sum_{uv \in E} \left( \rho_v^{uv}(x_{uv}) [v]+ \rho_u^{uv}(x_{uv}) [u]\right) \\
&=\sum_{uv \in E} \left(x_{uv}w_{vu}[v]+x_{uv}w_{uv}[u]\right) \\
&=\sum_{u \in V} \left(\sum_{v \in \Gamma_G(u)} x_{uv} w_{uv}\right) [u].
\end{split}
\end{equation*}
It follows that $\partial_1 c=0$ iff $\sum_{v \in \Gamma_G(u)} x_{uv} w_{uv}=0$ for all $u \in V$,
i.e. iff $(x_{uv})_{uv \in E} \in \calc$.
{\enp}

\end{document}